\documentclass[12pt,leqno]{article}
\usepackage[mathscr]{eucal}
\usepackage{amsmath,amssymb,latexsym,theorem,bbm,amsfonts}
\usepackage{verbatim}
\usepackage{color}
\usepackage{appendix}
\usepackage{array,epsfig,theorem,bbm,graphics,booktabs,url}

\newcommand{\NN}{\mathbb{N}}

\newcommand{\RR}{\mathbb{R}}

\newcommand{\ZZ}{\mathbb{Z}}

\newcommand{\tcY}{\widetilde{\cY}}

\newcommand{\cD}{{\mathcal D}}

\newcommand{\cY}{{\mathcal Y}}

\newcommand{\dd}{\mathrm{d}}
\newcommand{\ee}{\mathrm{e}}
\newcommand{\ff}{\mathrm{f}}
\newcommand{\ii}{\mathrm{i}}

\newcommand{\EE}{\operatorname{\mathbb{E}}}

\newcommand{\PP}{\operatorname{\mathbb{P}}}

\newcommand{\var}{\operatorname{Var}}
\newcommand{\cov}{\operatorname{Cov}}

\newcommand{\tS}{\widetilde{S}}

\newcommand{\vare}{\varepsilon}

\renewcommand{\mid}{\,|\,}

\renewcommand{\leq}{\leqslant}

\newcommand{\distr}{\stackrel{\cD}{\longrightarrow}}
\newcommand{\distrf}{\stackrel{\cD_\ff}{\longrightarrow}}

\newcommand{\bbone}{\mathbbm{1}}

\newcommand{\proofend}{\hfill\mbox{$\Box$}}

\numberwithin{equation}{section}

\theoremstyle{change} \theorembodyfont{\em}
\newtheorem{Lem}{Lemma.}[section]
\newtheorem{Thm}[Lem]{Theorem.}
\newtheorem{Pro}[Lem]{Proposition.}

\theorembodyfont{\rm}

\begin{document}

\begin{center}
 {\bfseries\Large
   Iterated scaling limits for aggregation of random coefficient AR(1) and INAR(1) processes}

\vskip0.5cm

 {\sc\large
  Fanni $\text{Ned\'enyi}^{*,\diamond}$,
  Gyula $\text{Pap}^{*}$}

\end{center}

\vskip0.2cm

\noindent
 * Bolyai Institute, University of Szeged,
     Aradi v\'ertan\'uk tere 1, H--6720 Szeged, Hungary.

\noindent e--mails:
                    nfanni@math.u-szeged.hu (F. Ned\'enyi),
                    papgy@math.u-szeged.hu (G. Pap).

\noindent $\diamond$ Corresponding author.

\vskip0.2cm


\renewcommand{\thefootnote}{}
\vspace*{0.2cm}

\vspace*{-10mm}

\begin{abstract}
We discuss joint temporal and contemporaneous aggregation of \ $N$ \ independent
 copies of strictly stationary AR(1) and INteger-valued AutoRegressive processes of order 1 (INAR(1)) with random coefficient \ $\alpha \in (0, 1)$ \ and idiosyncratic
 innovations.
Assuming that \ $\alpha$ \ has a density function of the form
 \ $\psi(x) (1 - x)^\beta$, \ $x \in (0, 1)$, \ with
 \ $\lim_{x\uparrow 1} \psi(x) = \psi_1 \in (0, \infty)$, \ different Brownian limit processes of
 appropriately centered and scaled aggregated partial sums are shown to exist in case \ $\beta=1$ \  when taking first the limit as \ $N \to \infty$ \ and then
 the time scale \ $n \to \infty$, \ or vice versa.
This paper completes the one of Pilipauskait\.e and Surgailis \cite{PilSur}, and Barczy, Ned\'enyi and Pap \cite{BarNedPap}, where the iterated limits are given for every other possible value of the parameter \ $\beta$ \ for the two types of models.
\end{abstract}

\section{Introduction}

The aggregation problem is concerned with the relationship between individual (micro)
 behavior and aggregate (macro) statistics.
There exist different types of aggregation.
The scheme of contemporaneous (also called cross-sectional) aggregation of
 random-coefficient AR(1) models was firstly proposed by Robinson \cite{Rob} and Granger
 \cite{Gra} in order to obtain the long memory phenomena in aggregated time series.

Puplinskait\.e and Surgailis \cite{PupSur1,PupSur2} discussed aggregation of
 random-coefficient AR(1) processes with infinite variance and innovations in the domain
 of attraction of a stable law.
Related problems for some network traffic models,
 \ $\mathrm{M}/\mathrm{G}/\infty$ \ queues with heavy-tailed activity periods, and
 renewal-reward processes have also been examined. 
On page 512 in Jirak \cite{Jir} one can find many references for papers dealing with the
 aggregation of continuous time stochastic processes, and the introduction of Barczy, Ned\'enyi and Pap \cite{BarNedPap}
contains a detailed overview on the topic.

The aim of the present paper is to complete the papers of Pilipauskait\.e and Surgailis
 \cite{PilSur} and Barczy, Ned\'enyi and Pap \cite{BarNedPap} by giving the appropriate iterated limit theorems for both the randomized AR(1) and INAR(1) models when the parameter \ $\beta=1$, \ which case is not investigated in both papers.

Let \ $\ZZ_+$, \ $\NN$, \ $\RR$ \ and \ $\RR_+$ \ denote the set of
 non-negative integers, positive integers, real numbers and non-negative real numbers, respectively.
The paper of Pilipauskait\.e and Surgailis
 \cite{PilSur} discusses the limit behavior of sums
 \begin{equation}\label{StNn}
  S_t^{(N,n)} := \sum_{j=1}^N \sum_{k=1}^{\lfloor nt \rfloor} X_k^{(j)} ,
  \qquad t \in \RR_+ ,
 \qquad N, n \in \NN ,
 \end{equation}
 where \ $(X_k^{(j)})_{k\in\ZZ_+}$, \ $j \in \NN$, \ are
 independent copies of a stationary random-coefficient AR(1) process
 \begin{equation}\label{RCAR}
  X_k = \alpha X_{k-1} + \vare_k , \qquad k \in \NN ,
 \end{equation}
 with standardized independent and identically distributed (i.i.d.) innovations
 \ $(\vare_k)_{k\in\NN}$ \ having \ $\EE(\vare_1) = 0$ \ and
 \ $\var(\vare_1) = 1$, \ and a random coefficient \ $\alpha$ \ with values in \ $[0, 1)$,
 \ being independent of \ $(\vare_k)_{k\in\NN}$ \ and admitting a probability
 density function of the form
 \begin{equation}\label{varphi}
  \psi(x) (1 - x)^\beta , \qquad x \in [0, 1) ,
 \end{equation}
 where \ $\beta \in (-1, \infty)$ \ and \ $\psi$ \ is an integrable
 function on \ $[0, 1)$ \ having
 a limit \ $\lim_{x\uparrow 1} \psi(x) = \psi_1 > 0$.
\ Here the distribution of \ $X_0$ \ is chosen as the unique stationary distribution of
 the model \eqref{RCAR}.
Its existence was shown in Puplinskait\.e and Surgailis \cite[Proposition 1]{PupSur1}.
We point out that they considered so-called idiosyncratic innovations, i.e., the
 innovations \ $(\vare^{(j)}_k)_{k\in \NN}$, \ $j \in \NN$, \ belonging  to
 \ $(X^{(j)}_k)_{k\in\ZZ_+}$, \ $j \in \NN$, \ are independent.
In Pilipauskait\.e and Surgailis \cite{PilSur} they derived scaling limits of the finite dimensional distributions of
 \ $(A_{N,n}^{-1} S_t^{(N,n)})_{t\in\RR_+}$, \ where \ $A_{N,n}$ \ are some scaling
 factors and first \ $N \to \infty$ \ and then \ $n \to \infty$, \ or vice versa, or
 both \ $N$ \ and \ $n$ \ increase to infinity, possibly with different rates.
 The iterated limit theorems for both orders of iteration are presented in the paper of  Pilipauskait\.e and Surgailis \cite{PilSur}, in Theorems 2.1 and 2.3, along with results concerning simultaneous limit theorems in Theorem 2.2 and 2.3. We note that the theorems cover different ranges of the possible values of \ $\beta\in (-1,\infty)$, \ namely, \ $\beta\in(-1,0)$, \ $\beta =0$, \   $\beta\in(0,1)$, and \ $\beta>1$. \ Among the limit processes is a fractional Brownian motion, lines with random slopes where the slope is a stable variable, a stable L\'evy process, and a Wiener process. Our paper deals with the missing case when \ $\beta=1$, \ for both two orders of iteration.

The paper of Barczy, Ned\'enyi and Pap \cite{BarNedPap} discusses the limit behavior of sums \eqref{StNn},
 where \ $(X_k^{(j)})_{k\in\ZZ_+}$, \ $j \in \NN$, \  are
 independent copies of a stationary random-coefficient INAR(1) process. The usual INAR(1) process with non-random-coefficient is defined as
 \begin{equation}\label{INAR1}
  X_k = \sum_{j=1}^{X_{k-1}} \xi_{k,j} + \vare_k, \qquad k \in \NN ,
 \end{equation}
 where \ $(\vare_k)_{k\in\NN}$ \ are i.i.d.\ non-negative integer-valued
 random variables, \ $(\xi_{k,j})_{k,j\in\NN}$ \ are i.i.d.\ Bernoulli random
 variables with mean \ $\alpha \in [0, 1]$, \ and \ $X_0$ \ is a non-negative
 integer-valued random variable such that \ $X_0$,
 \ $(\xi_{k,j})_{k,j\in\NN}$ \ and \ $(\vare_k)_{k\in\NN}$ \ are
 independent.
By using the binomial thinning operator \ $\alpha\,\circ$ \ due to Steutel and van Harn
 \cite{SteHar},
the INAR(1) model in \eqref{INAR1} can be considered as
 \begin{equation}\label{SteHarINAR1}
  X_k = \alpha \circ X_{k-1} + \vare_k , \qquad k \in \NN ,
 \end{equation}
which form captures the resemblance with the AR(1) model.
We note that an INAR(1) process can also be considered as a special branching process
 with immigration having Bernoulli offspring distribution.

We will consider a certain randomized INAR(1) process with randomized thinning parameter
 \ $\alpha$, \ given formally by the recursive equation \eqref{SteHarINAR1},
 where \ $\alpha$ \ is a random variable with values in \ $(0, 1)$.
\ This means that, conditionally on \ $\alpha$, \ the process \ $(X_k)_{k\in\ZZ_+}$ \ is an INAR(1)
 process with thinning parameter \ $\alpha$.
\ Conditionally on \ $\alpha$, \ the i.i.d.\ innovations
 \ $(\vare_k)_{k\in\NN}$ \ are supposed to have a Poisson distribution with
 parameter \ $\lambda \in (0, \infty)$, \ and the conditional distribution of the initial
 value \ $X_0$ \ given \ $\alpha$ \ is supposed to be the unique stationary distribution,
 namely, a Poisson distribution with parameter \ $\lambda/(1-\alpha)$.
\ For a rigorous construction of this process see Section 4 of Barczy, Ned\'enyi and Pap \cite{BarNedPap}. The iterated limit theorems for both orders of iteration ---that are analogous to the ones in case of the randomized AR(1) model--- are presented in the latter paper, in Theorems 4.6-4.12. This paper deals with the missing case when \ $\beta=1$, \ for both two orders of iteration.
\label{INAR}

\section{Iterated aggregation of randomized INAR(1) processes with
          Poisson innovations}
\label{RINAR}

Let \ $\alpha^{(j)}$, \ $j \in \NN$, \ be a sequence of independent copies of the
 random variable \ $\alpha$, \ and let \ $(X^{(j)}_k)_{k\in\ZZ_+}$, \ $j \in \NN$,
 \ be a sequence of independent copies of the process \ $(X_k)_{k\in\ZZ_+}$
 \ with idiosyncratic innovations (i.e., the innovations
 \ $(\vare^{(j)}_k)_{k\in\NN}$, $j\in\NN$, \ belonging to
 \ $(X^{(j)}_k)_{k\in\ZZ_+}$, \ $j \in \NN$, \ are independent) such that
 \ $(X^{(j)}_k)_{k\in\ZZ_+}$ \ conditionally on \ $\alpha^{(j)}$ \ is a strictly
 stationary INAR(1) process with Poisson innovations for all \ $j \in \NN$.

First we examine a simple aggregation procedure.
For each \ $N \in \NN$, \ consider the stochastic process
 \ $\tS^{(N)} = (\tS^{(N)}_k)_{k\in\ZZ_+}$ \ given by
 \[
   \tS^{(N)}_k
   := \sum_{j=1}^N \big(
		X^{(j)}_k - \EE(X^{(j)}_k \mid \alpha^{(j)})
		\big)
   = \sum_{j=1}^N \Bigl(X^{(j)}_k - \frac{\lambda}{1-\alpha^{(j)}}\Bigr) ,
   \qquad k \in \ZZ_+ .
 \]
The following two propositions are Proposition 4.1 and 4.2 of Barczy, Ned\'enyi and Pap \cite{BarNedPap}. 
We will use \ $\distrf$ \ or \ $\cD_\ff\text{-}\hspace*{-1mm}\lim$ \ for the weak
 convergence of the finite dimensional distributions.
\begin{Pro}\label{simple_aggregation_random}
If \ $\EE\bigl(\frac{1}{1-\alpha}\bigr) < \infty$, \ then
 \[
   N^{-\frac{1}{2}} \tS^{(N)}
   \distrf \tcY \qquad \text{as \ $N \to \infty$,}
 \]
 where \ $(\tcY_k)_{k\in\ZZ_+}$ \ is a stationary Gaussian process with zero mean
 and covariances
 \begin{equation}\label{covariance}
  \EE(\tcY_0 \tcY_k)
  = \cov\left(X_0 - \frac{\lambda}{1-\alpha},
              X_k - \frac{\lambda}{1-\alpha}\right)
  = \lambda \EE\Bigl(\frac{\alpha^k}{1-\alpha}\Bigr), \quad k\in\ZZ_+.
 \end{equation}
\end{Pro}


\begin{Pro}\label{simple_aggregation_random3}
We have
 \[
   \biggl(n^{-\frac{1}{2}}
          \sum_{k=1}^{\lfloor nt \rfloor}
           \tS^{(1)}_k\biggr)_{t\in\RR_+}
   = \biggl(n^{-\frac{1}{2}}
            \sum_{k=1}^{\lfloor nt \rfloor}
            (X^{(1)}_k - \EE(X^{(1)}_k \mid \alpha^{(1)}))\biggr)_{t\in\RR_+}
   \distrf \frac{\sqrt{\lambda(1+\alpha)}}{1-\alpha} B
 \]
 as \ $n \to \infty$, \ where \ $B = (B_t)_{t\in\RR_+}$ \ is a standard Brownian
 motion, independent of \ $\alpha$.
\end{Pro}

In the forthcoming theorems we assume that the
 distribution of the random variable \ $\alpha$, \ i.e., the mixing distribution, has a
 probability density described in \eqref{varphi}. We note that the form of this density function indicates \ $\beta>-1$.
\ Furthermore, if \ $\alpha$ \ has such a density function, then for each \ $\ell \in \NN$ \ the expectation
 \ $\EE ({(1-\alpha)^{-\ell}})$ \ is finite if and only if
 \ $\beta > \ell - 1$.

For each \ $N, n \in \NN$, \ consider the stochastic process
 \ $\widetilde S^{(N,n)} = (\widetilde S_t^{(N,n)})_{t\in\RR_+}$ \ given by
 \[
   \widetilde S_t^{(N,n)}
   := \sum_{j=1}^N \sum_{k=1}^{\lfloor nt \rfloor}
       \big(X^{(j)}_k - \EE(X^{(j)}_k \mid \alpha^{(j)})\big),
   \qquad t\in\RR_+.
 \]

\begin{Thm}\label{beta=1_1}
If \ $\beta=1$, \ then
 \[
   \cD_\ff\text{-}\hspace*{-1mm}\lim_{n\to\infty} \,
   \cD_\ff\text{-}\hspace*{-1mm}\lim_{N\to\infty} \,
    {(n\log n)}^{-\frac{1}{2}} N^{-\frac{1}{2}} \,
   \widetilde S^{(N,n)}
   = \sqrt{2\lambda \psi_1}B,
 \]
where \ $B=(B_t)_{t\in\RR_+}$ \ is a standard Wiener process.
\end{Thm}

\noindent{\bf Proof of Theorem \ref{beta=1_1}.}
Since \ $\EE((1-\alpha)^{-1})<\infty$,\ the condition in Proposition \ref{simple_aggregation_random} is satisfied, meaning that
\[
   N^{-\frac{1}{2}} \tS^{(N)}
   \distrf \tcY \qquad \text{as \ $N \to \infty$,}
 \]
 where \ $(\tcY_k)_{k\in\ZZ_+}$ \ is a stationary Gaussian process with zero mean
 and covariances
 \begin{equation*}
  \EE(\tcY_0 \tcY_k)
  = \cov\left(X_0 - \frac{\lambda}{1-\alpha},
              X_k - \frac{\lambda}{1-\alpha}\right)
  = \lambda \EE\Bigl(\frac{\alpha^k}{1-\alpha}\Bigr) , \qquad k \in \ZZ_+ .
 \end{equation*} 
Therefore, it suffices to show that
\[
   \cD_\ff\text{-}\hspace*{-1mm}\lim_{n\to\infty} \,
    \frac{1}{\sqrt{n\log n}}  \,
   \sum_{k=1}^{\lfloor nt \rfloor} \tcY_k
   = \sqrt{2\lambda \psi_1}B,
 \]
where \ $B=(B_t)_{t\in\RR_+}$ \ is a standard Wiener process.
This follows from the continuity theorem if for all \ $t_1,t_2\in \NN$ \ we have
\begin{equation}\label{covconv}
\cov\left(
   \frac{1}{\sqrt{n\log n}}  \,
   \sum_{k=1}^{\lfloor nt_1 \rfloor} \tcY_k,
   \frac{1}{\sqrt{n\log n}}  \,
   \sum_{k=1}^{\lfloor nt_2 \rfloor} \tcY_k
\right)\to 2\lambda \psi_1 \min(t_1,t_2),
\end{equation}
as \ $n\to \infty$.
\ By \eqref{covariance} we have
\begin{equation*}
\begin{split}
&
\cov\left(
   \frac{1}{\sqrt{n\log n}}  \,
   \sum_{k=1}^{\lfloor nt_1 \rfloor} \tcY_k,
   \frac{1}{\sqrt{n\log n}}  \,
   \sum_{k=1}^{\lfloor nt_2 \rfloor} \tcY_k
\right)
=\frac{\lambda}{n\log n}\EE\left(
\sum_{k=1}^{\lfloor nt_1\rfloor}
\sum_{\ell=1}^{\lfloor nt_2\rfloor}
\frac{\alpha^{|k-\ell|}}{1-\alpha}
\right)\\
&=
\frac{\lambda}{n\log n} \int_0^1 \sum_{k=1}^{\lfloor nt_1\rfloor}
\sum_{\ell=1}^{\lfloor nt_2\rfloor}
\frac{a^{|k-\ell|}}{1-a} \psi(a) (1-a) \, \dd a.
\end{split}
\end{equation*}
First we derive
\begin{equation}\label{conv}
\frac{1}{n\log n} \int_0^{1} \sum_{k=1}^{\lfloor nt_1\rfloor}
\sum_{\ell=1}^{\lfloor nt_2\rfloor}
a^{|k-\ell|} \, \dd a \to 2 \min(t_1,t_2),  
\end{equation}
as \ $n\to \infty$. \
Indeed, if we suppose that \ $t_2>t_1$, \ then
\begin{equation*}
\begin{split}
&\int_0^{1} \sum_{k=1}^{\lfloor nt_1\rfloor}
\sum_{\ell=1}^{\lfloor nt_2\rfloor}
a^{|k-\ell|} \, \dd a 
=\sum_{k=1}^{\lfloor nt_1\rfloor}
\sum_{\ell=1}^{\lfloor nt_2\rfloor}
\frac{1}{|k-\ell|+1}
\\
&=
(\lfloor nt_1\rfloor+1)
(H(\lfloor nt_1\rfloor)-1)+2-\lfloor nt_1\rfloor
+
\lfloor nt_1\rfloor
(H(\lfloor nt_2\rfloor)-1)\\
&
\quad +\big(\lfloor nt_2\rfloor-\lfloor nt_1\rfloor+1\big)\left(
H(\lfloor nt_2\rfloor)-H(\lfloor nt_2\rfloor-\lfloor nt_1\rfloor+1)\right)\\
&
= 
(\lfloor nt_1\rfloor+1)
(\log (\lfloor nt_1\rfloor)+O(1))+2-\lfloor nt_1\rfloor
+
\lfloor nt_1\rfloor
(\log \lfloor nt_2\rfloor +O(1))\\
&
\quad +\big(\lfloor nt_2\rfloor-\lfloor nt_1\rfloor+1\big)\left(
\log (\lfloor nt_2 \rfloor)-\log(\lfloor nt_2\rfloor-\lfloor nt_1\rfloor+1)+O(1)\right),
\end{split}
\end{equation*}
where \ $H(n)$ \ denotes the \ $n$\,-th harmonic number, and it is well known that \ $H(n)=\log n+O(1)$ \ for every \ $n\in \NN$. \ Therefore, convergence \eqref{conv} holds.
Consequently, \eqref{covconv} will follow from
\[
  I_n := \frac{1}{n\log n} \int_0^{1} \sum_{k=1}^{\lfloor nt_1\rfloor}
\sum_{\ell=1}^{\lfloor nt_2\rfloor}
a^{|k-\ell|} |\psi(a) - \psi_1| \, \dd a
\to 0
\]
as \ $n\to \infty$. \
Note that for every \ $\varepsilon>0$ \ there is a \ $\delta_\varepsilon>0$ \ such that for every \ $a\in(1-\delta_\varepsilon, 1)$ \ it holds that \ $|\psi(a)-\psi_1|<\varepsilon$.
\ Hence
\begin{align*}
 n \log n \,  I_n
 &\leq  \int_0^{1-\delta_\vare} \sum_{k=1}^{\lfloor nt_1\rfloor}
\sum_{\ell=1}^{\lfloor nt_2\rfloor}
a^{|k-\ell|} (\psi(a) + \psi_1) \, \dd a
+\int_{1-\delta_\vare}^{1} \sum_{k=1}^{\lfloor nt_1\rfloor}
\sum_{\ell=1}^{\lfloor nt_2\rfloor}
a^{|k-\ell|} |\psi(a) - \psi_1| \, \dd a \\
&\leq \int_0^{1-\delta_\vare}\frac{2\lfloor nt_1\rfloor}{\delta_\vare} (\psi(a) + \psi_1) \, \dd a
+
\vare \int_{1-\delta_\vare}^1\sum_{k=1}^{\lfloor nt_1\rfloor}
\sum_{\ell=1}^{\lfloor nt_2\rfloor}
a^{|k-\ell|} \, \dd a,
\end{align*}
meaning that for every \ $\vare>0$ \ by \eqref{conv} we have \ $ \limsup_{n\to \infty} |I_n|\leq 0+\varepsilon 4 \psi_1 \min(t_1,t_2)$, \ resulting that \ $  \lim_{n\to \infty}{I_n}=0$,\ which completes the proof.
\proofend

\begin{Thm}\label{beta=1_2}
If \ $\beta=1$, \ then 
 \[
   \cD_\ff\text{-}\hspace*{-1mm}\lim_{N\to\infty} \,
   \cD_\ff\text{-}\hspace*{-1mm}\lim_{n\to\infty} \,
    \frac{1}{\sqrt{n N\log N}}\,
   \widetilde S^{(N,n)}
   = \sqrt{\lambda \psi_1}B,
 \]
where \ $B=(B_t)_{t\in\RR_+}$ \ is a standard Wiener process.
\end{Thm}

\noindent{\bf Proof of Theorem \ref{beta=1_2}.}
By the second proof of Theorem 4.9 of Barczy, Ned\'enyi and Pap \cite{BarNedPap} it suffices to show that
\[
\frac{1}{N \log N}\sum_{j=1}^N\frac{\lambda(1+\alpha^{(j)})}{(1-\alpha^{( j)})^2} \distr \lambda \psi_1, \qquad N\to \infty.
\]
Let us apply Theorem 7.1 of Resnick \cite{Resnick} with
\[
X_{N, j}:=\frac{1}{N}\frac{\lambda(1+\alpha^{(j)})}{(1-\alpha^{(j)})^2},
\]
meaning that
\[
N\PP(X_{N,1}>x)=N\PP\left(\frac{\lambda(1+\alpha)}{(1-\alpha)^2}>Nx\right)=N\int_{1-\widetilde{h}(\lambda, Nx)}^{1} \psi(a)(1-a) \dd a,
\]
where \ $\widetilde{h}(\lambda, x)=(1/4+\sqrt{1/16+x/(2\lambda)})^{-1}.$ \
Note that for every \ $\varepsilon>0$ \ there is a \ $\delta_\varepsilon>0$ \ such that for every \ $a\in(1-\delta_\varepsilon, 1)$ \ it holds that \ $|\psi(a)-\psi_1|<\varepsilon$. Then, 
\[
N\int_{1-\widetilde{h}(\lambda, Nx)}^{1} |\psi(a)-\psi_1| (1-a) \dd a\leq
N \varepsilon \frac{(\widetilde{h}(\lambda, Nx))^2}{2}\leq \frac{\varepsilon \lambda}{x}
\] 
for every \ $x>0$ \ and large enough \ $N$.\ Therefore, for every \ $x>0$ \ we have
\begin{equation*}
\begin{split}
&
\lim_{N\to \infty} {N\PP(X_{N,1}>x)}=\lim_{N\to \infty} N\int_{1-\widetilde{h}(\lambda, Nx)}^{1} \psi_1(1-a) \dd a\\
&=\lim_{N\to \infty} 
N\psi_1\frac{(\widetilde{h}(\lambda, Nx))^2}{2}
=\lim_{N\to \infty} 
\frac{\psi_1}{2}\frac{N}{\left(\frac{1}{4}+\sqrt{\frac{1}{16}+\frac{Nx}{2\lambda}}\right)^2}=\frac{\psi_1\lambda}{x}=:\nu([x,\infty)),
\end{split}
\end{equation*}
where \ $\nu$ \ is obviously a L\'evy-measure.
By the decomposition
\begin{equation*}
N\EE\left(
X_{N,1}^2\bbone_{\{|X_{N,1}|\leq \varepsilon\}}
\right)=
N\int_0^{1-\widetilde{h}(\lambda, N\varepsilon)}\left(\frac{\lambda(1+a)}{N(1-a)^2}\right)^2\psi(a)(1-a)\dd a
= I_N^{(1)} + I_N^{(2)} ,
\end{equation*}
where
\[
I_N^{(1)}
:=
N\int_0^{1-\delta_\varepsilon}\left(\frac{\lambda(1+a)}{N(1-a)^2}\right)^2\psi(a)(1-a)\dd a\leq
\frac{1}{N}\lambda ^2\frac{2^2}{\delta_\varepsilon^4}1\to 0 
\]
as \ $N\to \infty$, \ and
\begin{equation*}
\begin{split}
I_N^{(2)}
&:=
N\int_{1-\delta_\varepsilon}^{1-\widetilde{h}(\lambda, N\varepsilon)}\left(\frac{\lambda(1+a)}{N(1-a)^2}\right)^2\psi(a)(1-a)\dd a
\\
&\leq 
\frac{8\psi_1\lambda^2}{N}  \int_{1-\delta_\varepsilon}^{1-\widetilde{h}(\lambda, N\varepsilon)}  \frac{\dd a}{(1-a)^3}
=
\frac{4\psi_1\lambda^2}{N}  \left[\widetilde{h}(\lambda, N\varepsilon)^{-2}
-\delta_\varepsilon^{-2}
\right]\leq 8\psi_1\lambda^2 \varepsilon
\end{split}
\end{equation*}
for large enough \ $N$ \ values, so
it follows that
\[
\lim_{\varepsilon\to 0}\limsup_{N\to \infty}N\EE\left(
X_{N,1}^2\bbone_{\{|X_{N,1}|\leq \varepsilon\}}
\right)=0.
\]
Therefore, by applying Theorem 7.1 of Resnick \cite{Resnick} with the choice \ $t=1$ \ we get that
\begin{equation*}
\begin{split}
&
\sum_{j=1}^N \left[\frac{\lambda (1+\alpha^{(j)})}{N(1-\alpha^{(j)})^2}-
\EE\left(
\frac{\lambda (1+\alpha)}{N(1-\alpha)^2}\bbone_{\left\{\frac{\lambda (1+\alpha)}{N(1-\alpha)^2}\leq 1\right\}}
\right)
\right]
\\
&
=
\sum_{j=1}^N \Bigg[\frac{\lambda (1+\alpha^{(j)})}{N(1-\alpha^{(j)})^2}-
\frac{\lambda\psi_1}{N}\int_0^{1-\sqrt{\frac{2\lambda}{N}}}\frac{2}{(1-a)^2}(1-a)\dd a
\\
&\quad +
\frac{\lambda\psi_1}{N}\int_0^{1-\sqrt{\frac{2\lambda}{N}}}\frac{2}{(1-a)^2}(1-a)\dd a-
\frac{\lambda\psi_1}{N}\int_0^{1-\widetilde{h}(\lambda,N)}\frac{2}{(1-a)^2}(1-a)\dd a
\\
&
\quad +
\frac{\lambda\psi_1}{N}\int_0^{1-\widetilde{h}(\lambda,N)}\frac{2}{(1-a)^2}(1-a)\dd a
-
\frac{\lambda\psi_1}{N}\int_0^{1-\widetilde{h}(\lambda,N)}\frac{1+a}{(1-a)^2}(1-a)\dd a
\\
&
\quad +
\frac{\lambda\psi_1}{N}\int_0^{1-\widetilde{h}(\lambda,N)}\frac{1+a}{(1-a)^2}(1-a)\dd a-
\frac{\lambda}{N}\int_0^{1-\widetilde{h}(\lambda,N)}\!\frac{1+a}{(1-a)^2}\psi(a)(1-a)\dd a
\Bigg]\\
&=:\frac{\lambda}{N}\sum_{j=1}^N J_{j,N}^{(0)}+\lambda J_N^{(1)}+\lambda J_N^{(2)}+\lambda J_N^{(3)}\distr X_0,
\end{split}
\end{equation*}
where by (5.37) of Resnick \cite{Resnick}
\[
\EE(\ee^{\ii \theta X_0})=\exp\left\{
\int_1^\infty (\ee^{\ii \theta x}-1)\frac{\psi_1 \lambda \dd x}{x^2}
+
\int_0^1(\ee^{\ii \theta x}-1-\ii\theta x) \frac{\psi_1\lambda \dd x}{x^2}
\right\}, \quad \theta\in\RR.
\]
We show that
\[
\frac{|J_N^{(1)}|+|J_N^{(2)}|+|J_N^{(3)}|}{\log N}\to 0, \qquad N\to \infty,
\]
resulting
\begin{equation*}
\begin{split}
&
\frac{1}{\log N}\sum_{j=1}^N\frac{\lambda (1+\alpha^{(j)})}{N(1-\alpha^{(j)})^2}=
\frac{1}{\log N}\sum_{j=1}^N
\left[\frac{\lambda (1+\alpha^{(j)})}{N(1-\alpha^{(j)})^2}-
\frac{\lambda\psi_1}{N}\int_0^{1-\sqrt{\frac{2\lambda}{N}}}\frac{2}{1-a}\dd a\right]\\
&
\qquad +
\frac{2\lambda \psi_1}{\log N}\left(-\log \left( \sqrt{\frac{2\lambda}{N}}\right)\right)
\distr 0 \cdot X_0+\lambda \psi_1=\lambda \psi_1, \qquad N\to \infty.
\end{split}
\end{equation*}
Indeed,
\[
\frac{J_N^{(1)}}{\log N}=\frac{\psi_1}{\log N} \int_{1-\sqrt{\frac{2\lambda}{N}}}^{1-\widetilde{h}(\lambda,N)} \frac{2}{1-a} \dd a=\frac{2\psi_1}{\log N} \log \left(
\sqrt{\frac{2\lambda}{N}}\left(
\frac{1}{4}+\sqrt{\frac{1}{16}+\frac{N}{2\lambda}}
\right)
\right)
\]
converges to 0 as \ $N\to \infty$. \ Moreover,
\[
\frac{J_N^{(2)}}{\log N}= \frac{\psi_1}{\log N}\int_0^{1-\widetilde{h}(\lambda,N)} \frac{1-a}{(1-a)^2}(1-a)\dd a = \frac{\psi_1}{\log N} \left(1-\frac{1}{\frac{1}{4}+\sqrt{\frac{1}{16}+\frac{N}{2\lambda}}}\right)
\]
converges to 0 as \ $N\to \infty$. \ 
Finally,
\begin{equation*}
\begin{split}
&
\left|\frac{J_N^{(3)}}{\log N}\right|=
\left|
\frac{1}{\log N} \int_0^{1-\widetilde{h}(\lambda,N)} \frac{1+a}{1-a} (\psi_1-\psi(a)) \dd a
\right|\\
&
\leq
\frac{1}{\log N} \int_0^{1-\delta_\varepsilon} \frac{2}{\delta_\varepsilon}(\psi_1+\psi(a))\dd a+\frac{1}{\log N} \int_{1-\delta_\varepsilon}^{1-\widetilde{h}(\lambda,N)} \frac{2}{1-a} \, \varepsilon \, \dd a
\\
&
\leq
\frac{1}{\log N}\frac{2}{\delta_\varepsilon}(\psi_1+\delta_\vare^{-1})+\frac{2\varepsilon}{\log N}\left[
\log \delta_\varepsilon+\log\left(\frac{1}{4}+\sqrt{\frac{1}{16}+\frac{N}{2\lambda}}\right).
\right],
\end{split}
\end{equation*}
One can easily see that for all \ $\varepsilon>0$,\ we get 
\ $ \limsup_{N\to \infty} |{J_N^{(3)}}/{\log N}|\leq 0+\varepsilon$, \ resulting that \ $  \lim_{N\to \infty}{J_N^{(3)}}/{\log N}=0$,\ which completes the proof.
\proofend

\section{Iterated aggregation of randomized AR(1) processes with Gaussian innovations}
\label{RAR}

Let \ $\alpha^{(j)}$, \ $j \in \NN$, \ be a sequence of independent copies of the
 random variable \ $\alpha$, \ and let \ $(X^{(j)}_k)_{k\in\ZZ_+}$, \ $j \in \NN$,
 \ be a sequence of independent copies of the process \ $(X_k)_{k\in\ZZ_+}$
 \ with idiosyncratic Gaussian innovations (i.e., the innovations
 \ $(\vare^{(j)}_k)_{k\in\ZZ_+}$, $j\in\NN$, \ belonging to
 \ $(X^{(j)}_k)_{k\in\ZZ_+}$, \ $j \in \NN$, \ are independent) having zero mean and variance \ $\sigma^2\in\RR_+$ \ such that
 \ $(X^{(j)}_k)_{k\in\ZZ_+}$ \ conditionally on \ $\alpha^{(j)}$ \ is a strictly
 stationary AR(1) process for all \ $j \in \NN$. A rigorous construction of this random-coefficient process can be given similarly as in case of the randomized INAR(1) process  detailed in Section 4 of Barczy, Ned\'enyi and Pap \cite{BarNedPap}.

First we examine a simple aggregation procedure.
For each \ $N \in \NN$, \ consider the stochastic process
 \ $\tS^{(N)} = (\tS^{(N)}_k)_{k\in\ZZ_+}$ \ given by
 \[
   \tS^{(N)}_k := \sum_{j=1}^N X^{(j)}_k , \qquad k \in \ZZ_+ .
 \]
The following two propositions are the counterparts of Proposition \ref{simple_aggregation_random} and \ref{simple_aggregation_random3}, and can be proven similarly as the two concerning the randomized INAR(1) process.
\begin{Pro}\label{ARsimple_aggregation_random}
If \ $\EE\bigl(\frac{1}{1-\alpha^2}\bigr) < \infty$, \ then
 \[
   N^{-\frac{1}{2}} \tS^{(N)}
   \distrf \tcY \qquad \text{as \ $N \to \infty$,}
 \]
 where \ $(\tcY_k)_{k\in\ZZ_+}$ \ is a stationary Gaussian process with zero mean
 and covariances
 \begin{equation*}
  \EE(\tcY_0 \tcY_k)
  = \cov(X_0, X_k)
  = \sigma^2 \EE\Bigl(\frac{\alpha^k}{1-\alpha^2}\Bigr) , \qquad k \in \ZZ_+ .
 \end{equation*}
\end{Pro}


\begin{Pro}\label{ARsimple_aggregation_random3}
We have
 \[
   \biggl(n^{-\frac{1}{2}}
          \sum_{k=1}^{\lfloor nt \rfloor}
           \tS^{(1)}_k\biggr)_{t\in\RR_+}
   = \biggl(n^{-\frac{1}{2}}
            \sum_{k=1}^{\lfloor nt \rfloor}
             X^{(1)}_k \biggr)_{t\in\RR_+}
   \distrf \frac{\sigma}{1-\alpha} B
 \]
 as \ $n \to \infty$, \ where \ $B = (B_t)_{t\in\RR_+}$ \ is a standard Brownian
 motion, independent of \ $\alpha$.
\end{Pro}

Again, we assume that the
 distribution of the random variable \ $\alpha$ \ has a
 probability density described in \eqref{varphi}. Note that for each \ $\ell \in \NN$ \ the expectation
 \ $\EE ({(1-\alpha^2)^{-\ell}})$ \ is finite if and only if
 \ $\beta > \ell - 1$.

For each \ $N, n \in \NN$, \ consider the stochastic process
 \ $\tS^{(N,n)} = (\tS_t^{(N,n)})_{t\in\RR_+}$ \ given by
 \[
   \widetilde S_t^{(N,n)} := \sum_{j=1}^N \sum_{k=1}^{\lfloor nt \rfloor} X^{(j)}_k ,
   \qquad t \in \RR_+ .
 \]

\begin{Thm}\label{ARbeta=1_1}
If \ $\beta=1$, \ then
 \[
   \cD_\ff\text{-}\hspace*{-1mm}\lim_{n\to\infty} \,
   \cD_\ff\text{-}\hspace*{-1mm}\lim_{N\to\infty} \,
    {(n\log n)}^{-\frac{1}{2}} N^{-\frac{1}{2}} \,
   \widetilde S^{(N,n)}
   =  \sqrt{\sigma^2 \psi_1} B,
 \]
where \ $B=(B_t)_{t\in\RR_+}$ \ is a standard Wiener process.
\end{Thm}

\noindent{\bf Proof of Theorem \ref{ARbeta=1_1}.}
Since \ $\EE((1-\alpha^2)^{-1})<\infty$,\ the condition in Proposition \ref{ARsimple_aggregation_random} is satisfied, meaning that
\[
   N^{-\frac{1}{2}} \tS^{(N)}
   \distrf \tcY \qquad \text{as \ $N \to \infty$,}
 \]
 where \ $(\tcY_k)_{k\in\ZZ_+}$ \ is a stationary Gaussian process with zero mean
 and covariances
 \begin{equation*}
  \EE(\tcY_0 \tcY_k)
  = \cov\left(X_0,
              X_k \right)
  = \sigma^2 \EE\Bigl(\frac{\alpha^k}{1-\alpha^2}\Bigr) , \qquad k \in \ZZ_+ .
 \end{equation*} Therefore, it suffices to show that
\[
   \cD_\ff\text{-}\hspace*{-1mm}\lim_{n\to\infty} \,
    \frac{1}{\sqrt{n\log n}}  \,
   \sum_{k=1}^{\lfloor nt \rfloor} \tcY_k
   = \sqrt{\sigma^2 \psi_1} B,
 \]
where \ $B=(B_t)_{t\in\RR_+}$ \ is a standard Wiener process.
This follows from the continuity theorem, if for all \ $t_1,t_2\in \NN$ \ we have
\[
\cov\left(
   \frac{1}{\sqrt{n\log n}}  \,
   \sum_{k=1}^{\lfloor nt_1 \rfloor} \tcY_k,
   \frac{1}{\sqrt{n\log n}}  \,
   \sum_{k=1}^{\lfloor nt_2 \rfloor} \tcY_k
\right)\to \sigma^2 \psi_1 \min(t_1,t_2), \qquad n\to \infty.
\]
It is known that
\begin{equation*}
\begin{split}
&
\cov\left(
   \frac{1}{\sqrt{n\log n}}  \,
   \sum_{k=1}^{\lfloor nt_1 \rfloor} \tcY_k,
   \frac{1}{\sqrt{n\log n}}  \,
   \sum_{k=1}^{\lfloor nt_2 \rfloor} \tcY_k
\right)
=\frac{\sigma^2}{n\log n}\EE\left(
\sum_{k=1}^{\lfloor nt_1\rfloor}
\sum_{\ell=1}^{\lfloor nt_2\rfloor}
\frac{\alpha^{|k-\ell|}}{1-\alpha^2}
\right)\\
&=
\frac{\sigma^2}{n\log n} \int_0^1 \sum_{k=1}^{\lfloor nt_1\rfloor}
\sum_{\ell=1}^{\lfloor nt_2\rfloor}
\frac{a^{|k-\ell|}}{1-a^2} \psi(a) (1-a) \dd a\\
&
=
\frac{\sigma^2}{n\log n} \int_0^1 \sum_{k=1}^{\lfloor nt_1\rfloor}
\sum_{\ell=1}^{\lfloor nt_2\rfloor}
a^{|k-\ell|} \psi(a) \dd a
-\frac{\sigma^2}{n\log n} \int_0^1 \sum_{k=1}^{\lfloor nt_1\rfloor}
\sum_{\ell=1}^{\lfloor nt_2\rfloor}
\frac{a^{|k-\ell|+1}}{1+a} \psi(a) \dd a
\end{split}
\end{equation*}
It was shown in the proof of Theorem \ref{beta=1_1} that
\[
\frac{\sigma^2}{n\log n} \int_0^1 \sum_{k=1}^{\lfloor nt_1\rfloor}
\sum_{\ell=1}^{\lfloor nt_2\rfloor}
a^{|k-\ell|} \psi(a) \dd a\to 2 \sigma^2 \psi_1 \min(t_1,t_2), \qquad n\to \infty.
\]
We are going to prove that
\[
\frac{\sigma^2}{n\log n} \int_0^1 \sum_{k=1}^{\lfloor nt_1\rfloor}
\sum_{\ell=1}^{\lfloor nt_2\rfloor}
\frac{a^{|k-\ell|+1}}{1+a} \psi(a) \dd a -
 \frac{\sigma^2}{n\log n} \int_0^1 \sum_{k=1}^{\lfloor nt_1\rfloor}
\sum_{\ell=1}^{\lfloor nt_2\rfloor}
\frac{a^{|k-\ell|}}{1+a} \psi(a) \dd a
\]
converges to 0 as \ $n\to \infty$, \
which proves our theorem.
Indeed, \ if \ $t_2>t_1$, \ then
\begin{equation*}
\begin{split}
&
\left| \sum_{k=1}^{\lfloor nt_1\rfloor} \sum_{\ell=1}^{\lfloor nt_2\rfloor}
\left(
\frac{a^{|k-\ell|+1}}{1+a} -\frac{a^{|k-\ell|}}{1+a}\right) \right|
=\frac{1}{1+a}
\left|\sum_{k=1}^{\lfloor nt_1\rfloor}
\left(a^k-(a+1)+a^{\lfloor nt_2\rfloor -k+1}\right)\right|
\\
&
=\frac{1}{1+a}
\left|
\frac{a(a^{\lfloor nt_1\rfloor} -1)}{a-1}-(a+1)\lfloor nt_1\rfloor+\frac{a^{\lfloor nt_2\rfloor+1}-a^{\lfloor nt_2\rfloor-\lfloor nt_1\rfloor+1}}{a-1}\right| \leq 4\lfloor nt_2\rfloor,
\end{split}
\end{equation*}
and as \ $\psi(a),\, a\in(0,1)$ \ is integrable,
\[
\frac{\sigma^2}{n\log n} \int_0^1 4\lfloor nt_2\rfloor\psi(a) \dd a \to 0, \qquad n\to \infty.
\]
This completes the proof.
\proofend

\begin{Thm}\label{ARbeta=1_2}
If \ $\beta=1$, \ then 
 \[
   \cD_\ff\text{-}\hspace*{-1mm}\lim_{N\to\infty} \,
   \cD_\ff\text{-}\hspace*{-1mm}\lim_{n\to\infty} \,
    \frac{1}{\sqrt{n N\log N}}\,
   \widetilde S^{(N,n)}
   = \sqrt{\frac{\sigma^2 \psi_1}{2}} B,
 \]
where \ $B=(B_t)_{t\in\RR_+}$ \ is a standard Wiener process.
\end{Thm}

The proof is similar to the INAR(1) case since the only difference is a missing \ $1+\alpha$ \ factor in the numerator and the constants.

\bibliographystyle{plain}
\bibliography{./aggr2}
\end{document}